\documentclass[11pt,reqno]{amsart}
\usepackage{latexsym,amscd}
\usepackage{amsmath,amssymb,amsthm,url,mathrsfs}
\usepackage{enumerate,stmaryrd,scalefnt}
\usepackage{color}

\usepackage[colorlinks=true,urlcolor=black,linkcolor=black,citecolor=black]{hyperref}

\theoremstyle{plain}

\newtheorem{lemma}{Lemma}

\theoremstyle{definition}

\newcounter{claim}

\newcounter{conjecture}

 \theoremstyle{remark}

\newcommand{\<}{\ensuremath{\langle}}
\renewcommand{\>}{\ensuremath{\rangle}}
\renewcommand{\phi}{\ensuremath{\varphi}}

\newcommand{\lb}{\ensuremath{\llbracket}}
\newcommand{\rb}{\ensuremath{\rrbracket}}

\newcommand{\bA}{\ensuremath{\mathbf{A}}}
\newcommand{\bB}{\ensuremath{\mathbf{B}}}
\newcommand{\bC}{\ensuremath{\mathbf{C}}}

\renewcommand{\leq}{\ensuremath{\leqslant}}

\renewcommand{\lneq}{\ensuremath{\lneqq}}

\newcommand{\Con}{\ensuremath{\operatorname{Con}}}
\newcommand{\Sub}{\ensuremath{\operatorname{Sub}}}
\newcommand{\NSub}{\ensuremath{\operatorname{NSub}}}

\renewcommand{\phi}{\ensuremath{\varphi}}

\begin{document}
\title{Isotopic algebras with non-isomorphic congruence lattices}
\date{February 9, 2013}
\author[W. DeMeo]{William DeMeo}
\email{williamdemeo@gmail.com}
\urladdr{\url{http://williamdemeo.wordpress.com}}
\address{Department of Mathematics\\
University of South Carolina\\Columbia\\ USA}

\keywords{isotopy, isotopic algebras, congruence lattice}

\begin{abstract}
We give examples of pairs of isotopic algebras
with non-isomorphic congruence lattices. This answers
negatively the question of whether all isotopic algebras have isomorphic
congruence lattices.
\end{abstract}

\maketitle
\setcounter{section}{+1}
\setcounter{secnumdepth}{2}

It is well known that two algebras in a congruence modular variety that are
isotopic have isomorphic congruence lattices.  In fact, this result holds more
generally, but to date the most general result of this kind (recalled below)
assumes some form of congruence modularity. 
It is natural to ask to what extent the congruence modularity hypothesis could
be relaxed and whether it is possible to prove that all isotopic algebras have
isomorphic congruence lattices.   In this note we show
that the full generalization is not possible; we construct a class of
counter-examples involving pairs of algebras whose
congruence lattices are obviously not isomorphic, and then prove that these
pairs of algebras are isotopic.

If $\bA$, $\bB$, and $\bC$ are algebras of the same signature, we say that
$\bA$ and $\bB$ are \emph{isotopic over} $\bC$, and we write $\bA\sim_{\bC}\bB$,
if there exists an isomorphism $\phi: \bA \times \bC \rightarrow \bB \times \bC$ such that 
for all $a\in A,\, c\in C$, the second coordinate of $\phi(a,c)$ is $c$; that
is, $\phi(a,c) = (b,c)$ for some $b\in B$.
We say that $\bA$ and $\bB$ are \emph{isotopic}, and we write $\bA\sim \bB$, provided
$\bA\sim_{\bC}\bB$ for some $\bC$.  It is not hard to check that $\sim$ is an
equivalence relation.

If $\bA\sim_{\bC}\bB$ and the congruence lattice of $\bA \times \bC$ happens to
be modular, then we write $\bA \sim^{\mathrm{mod}}_{\bC} \bB$, in which case we say that
$\bA$ and $\bB$ are \emph{modular isotopic over} $\bC$.
We call $\bA$ and $\bB$ \emph{modular isotopic in one step},  denoted 
$\bA \sim^{\mathrm{mod}}_1 \bB$,
if they are modular isotopic over $\bC$ for some $\bC$. Finally, 
$\bA$ and $\bB$ are \emph{modular isotopic}, denoted 
$\bA \sim^{\mathrm{mod}} \bB$ if the pair $(\bA, \bB)$ belongs to the transitive
closure of $\sim^{\mathrm{mod}}_1$.

Let $\Con \bA$ denote the congruence lattice of $\bA$.  It is well known that
$\bA \sim^{\mathrm{mod}} \bB$ implies $\Con \bA \cong \Con \bB$.  The proof
of this result appearing in~\cite{alvi:1987} 
is a straight forward application of Dedekind's Transposition Principle.  
Since a version of this principle has been shown to 
hold even in the non-modular case (\cite{DTP}), we might hope that the
proof technique used in \cite{alvi:1987} could be used to show that  
$\bA\sim_{\bC}\bB$ implies $\Con \bA \cong \Con \bB$.  But this strategy quickly
breaks down, and the application of the perspectivity map, which works fine when
$\Con (\bA\times \bC)$ is modular, can fail if
$\Con (\bA\times \bC)$ is non-modular, even in cases where $\bA \cong \bB$. 

This note describes a class of examples in which $\bA\sim_{\bC}\bB$ and $\Con \bA
\ncong \Con \bB$. Although we suspect simpler examples can be found, the
construction given here is not complicated and reveals that congruence lattices of
isotopic algebras can differ quite dramatically.

For any group $G$, let $\Sub(G)$ denote the lattice of subgroups of $G$, and
let $\NSub(G)$ denote the lattice of normal subgroups of $G$.
A group $G$ is called a \emph{Dedekind group} provided every subgroup of $G$ is
normal ---i.e., $\Sub(G) = \NSub(G)$.  
We call $G$ a \emph{non-Dedekind group} if it has a non-normal subgroup; of
course this requires $G$ be nonabelian, but that is not
sufficient.  For example, the eight element quaternion group 
is a Dedekind group.

Let $S$ be any group and let $D$ denote the \emph{diagonal subgroup} of 
$S\times S$; that is, $D = \{(x,x) \mid x\in S\}$.
The \emph{filter above $D$ in} $\Sub(S\times S)$, which we denote by
$\lb D, S\times S\rb$, 
consists of the subgroups of $S\times S$ that contain $D$.  In symbols,
$\lb D, S\times S\rb = \{K \mid D \leq K \leq S\times S\}$.  
This is a sublattice of $\Sub(S\times S)$ and is described by the following easy lemma:
\begin{lemma}
\label{lem:1}
  The filter above the diagonal subgroup in the subgroup lattice of $S\times S$
  is isomorphic to the lattice of normal subgroups of $S$. In symbols, 
$\lb D, S\times S\rb \cong \NSub(S)$.
\end{lemma}

\noindent {\bf The example.}
Let $S$ be any finite non-Dedekind group. Let $G = S_1 \times S_2$, where
$S_1 \cong S_2 \cong S$, and let $D = \{(x_1,x_2)\in G \mid x_1 = x_2\}$,  the
diagonal subgroup of $G$. For ease of notation, put $T_1 = S_1 \times \{1\}$ and $T_2 = \{1\}\times S_2$.
Then $D \cong T_1 \cong T_2$, and these three 
subgroups are pair-wise compliments: $T_1\cap D = D\cap T_2 = T_1 \cap T_2 =
\{(1,1)\}$ and $ \<T_1, T_2\> = \<T_1,D\> = \<D, T_2\> = G$.

Let $\bA$ be the algebra whose universe is the set $A = G/T_1$ of left
cosets of $T_1$ in $G$, and whose operations are left multiplication by elements
of $G$. 
That is $\bA = \< G/T_1, G^{\bA}\>$ 
where, for each $g\in G$, the operation $g^{\bA} \in G^{\bA}$ is defined by
$g^{\bA}(xT_1) = (gx)T_1$.
Define the algebra $\bC = \< G/T_2, G^{\bC}\>$ similarly.  

The algebra $\bB$
will have as its universe the set $B = G/D$ of left cosets of $D$ in $G$, but in
this case we define the action of $G$ on $B$ with a slight twist:  for
each $g = (g_1, g_2) \in G$,
for each $(x_1, x_2)D \in G/D$,
\[
g^{\bB}((x_1,x_2)D) =  (g_2x_1, g_1 x_2)D.
\]
Let $\bB = \< G/D, G^{\bB}\>$, where $G^{\bB} =  \{g^{\bB} \mid g\in G\}$.

Consider the binary relation 
$\phi \subseteq (\bA \times \bC) \times (\bB \times \bC)$ defined by associating
to each ordered pair $((x_1,x_2)T_1, (y_1,y_2)T_2) \in \bA \times \bC$
 the pair $((x_2, y_1)D, (y_1,y_2)T_2) \in \bB \times \bC$.
Our claim (proved below) is that this relation defines a function $\phi \colon \bA \times \bC \rightarrow \bB
\times \bC$, and that this function is an isomorphism.
Since the second coordinates of $\phi$-related pairs are the same, this will
establish that $\bA\sim_{\bC}\bB$.  

Note that $\Con \bA$ is isomorphic to
the filter above $T_1$ in the subgroup lattice of $G$ (see~\cite[Lemma
4.20]{alvi:1987}), and this filter is isomorphic to 
the subgroup lattice of $S$. 
Thus, $\Con \bA\cong \Sub(S)$.
On the other hand, by
Lemma~\ref{lem:1} we have $\Con \bB \cong \NSub(S)$, and since we chose $S$ to
be a non-Dedekind group, we can conclude that 
$\Con \bB \cong\NSub(S) \lneq\Sub(S) \cong  \Con \bA$. 

The foregoing describes a class of examples, indexed by the group $S$.  The
group $S$ must be a non-Dedekind group but is otherwise arbitrary and can be
chosen so that $\Con\bA$ and $\Con\bB$ are not only non-isomorphic, but also 
very different in size. For example, if $S$ is a finite nonabelian simple
group, then $\Con\bB$ has just two elements, while
$\Con \bA \cong \Sub(S)$ can be enormous.

We conclude with the easy proofs of three claims which establish that $\phi$ is an isomorphism.
\\[4pt]
{\it Claim 1:} $\phi$ is a function.\\[2pt]
{\it Proof:} Suppose $(xT_1, yT_2) = (x'T_1, y'T_2)$.  Note that 
$(x_1,x_2)T_1 = (x_1',x_2')T_1$ and $T_1 = S_1\times \{1\}$ imply $x_2 =
x_2'$.  Similarly, 
$(y_1,y_2)T_2 = (y_1',y_2')T_2$ and $T_2 = \{1\}\times S_2$ imply $y_1 =
y_1'$. 
This yields
$((x_2,y_1)D, (y_1,y_2)T_2)
= ((x'_2,y'_1)D, (y'_1,y'_2)T_2)$, which proves that if $(a_1,b_1)\in \phi$ and
$(a_2,b_2) \in \phi$ and $a_1 = a_2$, then $b_1 = b_2$, so $\phi$ is a function.\\[4pt]
{\it Claim 2:} $\phi$ is a homomorphism.\\[2pt]
{\it Proof:} For
$g = (g_1, g_2) \in G$, for
$xT_1 = (x_1, x_2)T_1 \in G/T_1$, and for
$yT_2 = (y_1, y_2)T_2 \in G/T_2$, we have
\begin{align*}
\phi(g^{\bA\times \bC}(xT_1, yT_2)) 
&= \phi((g_1x_1,g_2x_2)T_1, (g_1y_1,g_2y_2)T_2) \\
&= ((g_2x_2,g_1y_1)D, (g_1y_1,g_2y_2)T_2) \\
&= (g^{\bB}((x_2,y_1)D), g^{\bC}(yT_2)) \\
&= g^{\bB\times \bC}((x_2,y_1)D, yT_2) \\
&= g^{\bB\times \bC}\phi(xT_1, yT_2).
\end{align*}
{\it Claim 3:} $\phi$ is bijective.\\[2pt]
{\it Proof:} 
Since $\phi$ is a function from the finite set $A$ to the finite set $B$, and
since $|A|= |G:T_1| = |G:D| = |B|$, it suffices to prove that $\phi$ is
injective.  Suppose $\phi(xT_1, yT_2) = \phi(x'T_1, y'T_2)$.
That is,  suppose
\[
((x_2,y_1)D, yT_2) = 
((x'_2,y'_1)D, y'T_2).
\]
Then, as above,
$(y_1,y_2)T_2 = (y'_1,y'_2)T_2$ implies $y_1 = y_1'$.  Therefore,
 since $(x_2,y_1)^{-1}(x'_2,y'_1) \in D$, we have
$x_2^{-1}x'_2 = y_1^{-1}y'_1 = 1$,  
so 
$(x_1,x_2)$ and $(x'_1,x'_2)$ can 
only differ in the first coordinate.  It follows that
$(x_1,x_2)T_1 = (x'_1,x_2)T_1$.
\def\cprime{$'$} \def\cprime{$'$}
  \def\ocirc#1{\ifmmode\setbox0=\hbox{$#1$}\dimen0=\ht0 \advance\dimen0
  by1pt\rlap{\hbox to\wd0{\hss\raise\dimen0
  \hbox{\hskip.2em$\scriptscriptstyle\circ$}\hss}}#1\else {\accent"17 #1}\fi}

\end{document}